\documentclass[10pt,twocolumn]{article}
\usepackage{amsmath, amsfonts, amssymb, amsthm}
\usepackage{times}
\usepackage{latex8}
\usepackage{fancyhdr}

\newtheorem{lemma}{Lemma}
\newtheorem{theorem}[lemma]{Theorem}
\newtheorem{corollary}[lemma]{Corollary}
\newtheorem{definition}[lemma]{Definition}
\newtheorem{proposition}[lemma]{Proposition}

\newtheorem{conjecture}[lemma]{Conjecture}
\numberwithin{lemma}{section}
\numberwithin{equation}{section}
\numberwithin{figure}{section}
\newcommand\C{{\mathbb{C}}}

\newcommand\Z{{\mathbb{Z}}}

\newcommand\set[1]{\{#1\}}
\newcommand\Sym{\operatorname{\textup{Sym}}}
\newcommand\Ind{\operatorname{\textup{Ind}}}
\newcommand\Cyc{\textup{Cyc}}
\newcommand\Symn{\textup{Sym}_n}
\newcommand{\Chart}{{\mathcal{C}}}
\newcommand{\HypGph}{{\mathcal{H}}}

\newcommand{\remove}[1]{{}}

\begin{document}

\title{Group-theoretic Algorithms for Matrix Multiplication}

\author{Henry Cohn\thanks{Microsoft Research, One Microsoft Way, Redmond,
WA 98052-6399, cohn@microsoft.com} \and Robert
Kleinberg\thanks{Department of Mathematics, Massachusetts
Institute of Technology, Cambridge, MA 02139, and Department of
Computer Science, Cornell University, Ithaca, NY 14853,
rdk@cs.cornell.edu.} \and Bal\'azs Szegedy\thanks{Microsoft
Research, One Microsoft Way, Redmond, WA 98052-6399,
szegedy@microsoft.com} \and Christopher Umans\thanks{Department of
Computer Science, California Institute of Technology, Pasadena, CA
91125, umans@cs.caltech.edu; supported in part by NSF grant
CCF-0346991 and an Alfred P. Sloan Research Fellowship.}}

\maketitle
\thispagestyle{empty}

\begin{abstract}
We further develop the group-theoretic approach to fast matrix
multiplication introduced by Cohn and Umans, and for the first
time use it to derive algorithms asymptotically faster than the
standard algorithm. We describe several families of wreath product
groups that achieve matrix multiplication exponent less than $3$,
the asymptotically fastest of which achieves exponent $2.41$. We
present two conjectures regarding specific improvements, one
combinatorial and the other algebraic. Either one would imply that
the exponent of matrix multiplication is $2$.
\end{abstract}

\vspace{-.1094132in}

\Section{Introduction}

\thispagestyle{fancy}

\fancyfoot[LE,LO]{\ \\ \parbox{6.875in}{\scriptsize{\ \\ Copyright
\copyright\ 2005 IEEE. Reprinted from Proceedings of the 46th
Annual Symposium on Foundations of Computer Science. This material
is posted here with permission of the IEEE.  Such permission of
the IEEE does not in any way imply IEEE endorsement of any of
Cornell University's products or services.  Internal or personal
use of this material is permitted.  However, permission to
reprint/republish this material for advertising or promotional
purposes or for creating new collective works for resale or
redistribution must be obtained from the IEEE by writing to
pubs-permissions@ieee.org.  By choosing to view this document, you
agree to all provisions of the copyright laws protecting it.}}}

The task of multiplying matrices is one of the most fundamental
problems in algorithmic linear algebra. Matrix multiplication
itself is a important operation, and its importance is magnified
by the number of similar problems that are reducible to it.

Following Strassen's discovery \cite{S} of an algorithm for $n
\times n$ matrix multiplication in $O(n^{2.81})$ operations, a
sequence of improvements has achieved ever better bounds on the
{\em exponent of matrix multiplication}, which is the smallest
real number $\omega$ for which $n \times n$ matrix multiplication
can be performed in $O(n^{\omega + \varepsilon})$ operations for
each $\varepsilon > 0$. The asymptotically fastest algorithm known
is due to Coppersmith and Winograd \cite{CW}, and it proves that
$\omega < 2.376$. Since 1990, there have been no better upper
bounds proved on $\omega$, although it is widely believed that
$\omega = 2$.

Recently, Cohn and Umans \cite{CU} proposed a new group-theoretic
approach to devising matrix multiplication algorithms. In this
framework, one selects a finite group $G$ satisfying a certain
property that allows $n \times n$ matrix multiplication to be
reduced to multiplication of elements of the {\em group algebra}
$\C[G]$. This latter multiplication is performed via a Fourier
transform, which reduces it to several smaller matrix
multiplications, whose sizes are the \textit{character degrees\/}
of $G$. This naturally gives rise to a recursive algorithm whose
running time depends on the character degrees. Thus the problem of
devising matrix multiplication algorithms in this framework is
imported into the domain of group theory and representation
theory.

One of the main contributions of \cite{CU} was to demonstrate that
several diverse families of non-abelian groups support the
reduction of $n \times n$ matrix multiplication to group algebra
multiplication. These include, in particular, families of groups
of size $n^{2 + o(1)}$. The existence of such families is a
necessary condition for the group-theoretic approach to prove
$\omega=2$, although it is not sufficient.

The main question raised in \cite{CU} is whether the proposed
approach could prove nontrivial bounds on $\omega$, i.e., prove
$\omega < 3$. This was shown to be equivalent to a question in
representation theory, Question 4.1 in \cite{CU}: is there a group
$G$ with subsets $S_1, S_2, S_3$ that satisfy the {\em triple
product property} (see Definition~\ref{definition:realize} below),
and for which $|S_1||S_2||S_3|
> \sum_i{d_i^3}$, where $\{d_i\}$ is the set of character degrees of $G$?

In this paper we resolve this question in the affirmative, which
immediately gives a simple matrix multiplication algorithm in the
group-theoretic framework that has running time $O(n^{2.9088})$.
The group we construct for this purpose is a {\em wreath product},
and in subsequent sections we describe similar constructions that
produce algorithms with running times $O(n^{2.48})$ and
$O(n^{2.41})$.

The main challenge in each case is to describe the three subsets
of the group that satisfy the triple product property. We give two
ways of organizing these descriptions, both of which give rise to
the $O(n^{2.48})$ algorithm relatively simply. We also advance two
natural conjectures related to these formulations, each of which
would imply that $\omega = 2$. The first is a combinatorial
conjecture (Conjecture~\ref{conjecture:strongUSP}), and the second
is an algebraic conjecture (Conjecture~\ref{conjecture:sdpp}).

The three subsets underlying the $O(n^{2.41})$ algorithm are
described in terms of a combinatorial object we call a {\em
Uniquely Solvable Puzzle} (or USP), which is a weakening of the
combinatorial object in our first conjecture.  An {\em optimal}
USP construction can be extracted from Coppersmith and Winograd's
paper \cite{CW}.

\fancyfoot[CE,CO]{\thepage}
\fancyfoot[LE,LO]{}

In fact, the reader familiar with Strassen's 1987 paper \cite{S87}
and Coppersmith and Winograd's paper \cite{CW} (or the
presentation of this material in, for example, \cite{BCS}) will
recognize that our exponent bounds of $2.48$ and $2.41$ match
bounds derived in those works. It turns out that with some effort
the algorithms in \cite{S87} and \cite{CW}, including Coppersmith
and Winograd's $O(n^{2.376})$ algorithm, all have analogues in our
group-theoretic framework.  The translation does not appear to be
systematic: the algorithms are based on similar principles, but in
fact they are not identical (the actual operations performed on
matrix entries do not directly correspond), and we know of no
group-theoretic interpretations of any earlier algorithms.  We
defer a complete account of this connection to the full version of
this paper.

We believe that, compared to existing algorithms, our
group-theoretic algorithms are simpler to state and simpler to
analyze.  They are situated in a clearer conceptual and
mathematical framework, in which, for example, the two conjectures
mentioned above are natural and easy to identify. Finally, they
avoid various complications of earlier algorithms.  For example,
they substitute the discrete Fourier transform, together with some
elementary facts in representation theory, for the seemingly ad
hoc trilinear form identities in introduced in \cite{S87}, and
they completely avoid the need to deal with degenerations and
border rank of tensors.

\paragraph{Outline.}
In the rest of this section, we establish notation and review
background from \cite{CU} on the group-theoretic approach to fast
matrix multiplication.  Section~\ref{section:cubesum} describes
the simplest group we have found that can prove a nontrivial bound
on the exponent of matrix multiplication.  In
Sections~\ref{section:usp} and~\ref{section:sdpp}, we carry out a
more elaborate construction in two different ways, each of which
has the potential of reaching $\omega=2$. The most fundamental
conceptual contribution in this paper is the \textit{simultaneous
triple product property\/}, which we introduce in
Section~\ref{section:stpp}.  It extends the triple product
property from \cite{CU}, and it encompasses and illuminates all of
our other constructions, as we explain in
Section~\ref{section:usingstpp}. Finally, in
Section~\ref{section:wreath} we show that any bound provable via
the simultaneous triple product property can in fact be proved
using only the approach of \cite{CU}.

\SubSection{Preliminaries and notation}

As usual $\omega$ denotes the exponent of matrix multiplication
over $\C$.

The set $\{1,2,\dots,k\}$ is denoted $[k]$. We write $ A \setminus
B = \{a \in A : a \not\in B\} $ and if $A$ and $B$ are subsets of
an abelian group we set $ A - B = \{a-b : a \in A, b \in B\}. $

The cyclic group of order $k$ is denoted $\Cyc_k$ (with additive
notation for the group law), and the symmetric group on a set $S$
is denoted $\Sym(S)$ (or $\Symn$ instead of $\Sym([n])$). If $G$
is a group and $R$ is a ring, then $R[G]$ will denote the group
algebra of $G$ with coefficients in $R$.

\remove{ When we discuss a group action of $G$ on $H$, it will
always be a left action (with the action of $g$ on $h$ written $g
\cdot h$) unless otherwise specified. For a right action of $G$ on
$H$, we write $h^g$ for the action of $g$ on $h$. }

When we discuss a group action, it will always be a left action
unless otherwise specified. If $G$ and $H$ are groups with a left
action of $G$ on $H$ (where the action of $g$ on $h$ is written $g
\cdot h$), then the semidirect product $H \rtimes G$ is the set $H
\times G$ with the multiplication law $ (h_1,g_1),(h_2,g_2) =
(h_1(g_1 \cdot h_2),g_1g_2). $ We almost always identify $H$ with
the subset $H \times \{1\}$ and $G$ with $\{1\} \times G$, so that
$(h,g)$ simply becomes the product $hg$. For a right action of $G$
on $H$, with the action of $g$ on $h$ written $h^g$, the
semidirect product $G \ltimes H$ is $G \times H$ with the
multiplication law $ (g_1,h_1)(g_2,h_2) = (g_1g_2, h_1^{g_2} h_2).
$ As in the previous case we identify $G$ and $H$ with the
corresponding subsets of $G \ltimes H$.

Other than for Lemma~\ref{lemma:wreath-char-degrees}, which is not
required for the main results of this paper, we will use only the
following basic facts from representation theory.  The group
algebra $\C[G]$ of a finite group $G$ decomposes as the direct
product $ \C[G] \cong  \C^{d_1 \times d_1} \times \dots \times
\C^{d_k \times d_k} $ of matrix algebras of orders
$d_1,\dots,d_k$. These orders are the character degrees of $G$, or
the dimensions of the irreducible representations.  It follows
from computing the dimensions of both sides that $|G| = \sum_i
d_i^2$.  It is also easy to prove that if $G$ has an abelian
subgroup $A$, then all the character degrees of $G$ are less than
or equal to the index $[G:A]$ (Proposition~2.6 in \cite{H}) . See
\cite{JL} and \cite{H} for further background on representation
theory.

The following elementary lemma (proof omitted) will prove useful
several times:

\begin{lemma}
\label{lemma:geom->arith} Let $s_1, s_2, \dots, s_n$ be
nonnegative real numbers, and suppose that for every vector $\mu =
(\mu_1,\dots,\mu_n)$ of nonnegative integers for which
$\sum_{i=1}^n \mu_i = N$ we have $ \binom{N}{\mu} \prod_{i = 1}^n
s_i^{\mu_i} \le C^N. $ Then $\sum_{i=1}^n s_i \le C$.
\end{lemma}

\remove{
\begin{proof}
For each probability distribution $p = (p_1,\dots,p_n)$, we can
let $N$ tend to infinity and choose $\mu$ so that $\lim_{N
\rightarrow \infty} \mu/N = p$.  As $N \to \infty$, the inequality
in the hypothesis of the lemma yields
$$
-\sum_i p_i \log p_i + \sum_i p_i \log s_i \le \log C
$$
after taking the $N$-th root and the logarithm. Setting $S =
\sum_i s_i$ and $p_i = s_i/S$ proves $\log S \le \log C$, as
desired.
\end{proof}
}

Occasionally we will need to bound the character degrees of wreath
products:

\begin{lemma}
\label{lemma:wreath-char-degrees} Let $\{d_k\}$ be the character
degrees of a finite group $H$ and let $\{c_j\}$ be the character
degrees of $\Symn \ltimes H^n$ (where $\Symn$ acts by permuting
the coordinates). Then $\sum_j c_j^\omega \le (n!)^{\omega - 1}
\left (\sum_k d_k^\omega \right )^n.$
\end{lemma}

\begin{proof}[Sketch of proof.]
When $H$ is abelian, the theorem follows from the elementary facts
that the character degrees of $\Symn \ltimes H^n$ are at most $n!$
(which is the index of $H^n$ in $\Symn \ltimes H^n$) and that
$\sum_j c_j^2 = |\Symn \ltimes H^n|$. For general $H$, we the
theorem can be derived from well-known characterizations of the
character degrees of $\Symn \ltimes H^n$ (see, e.g., Theorem~25.6
in \cite{H}).
\end{proof}

\remove{
\begin{proof}
When $H$ is abelian, we can use the elementary facts that the
character degrees of $\Symn \ltimes H^n$ are at most $n!$ (which
is the index of $H^n$ in $\Symn \ltimes H^n$) and that $\sum_j
c_j^2 = |\Symn \ltimes H^n|$ to obtain
$$
\sum_j c_j^\omega \le (n!)^{\omega - 2}\sum_j c_j^2 = (n!)^{\omega
- 1} |H|^n.
$$

For general $H$, we need further information regarding the
character degrees of $\Symn \ltimes H^n$. From Theorem~25.6 in
\cite{H} we get the following description: The symmetric group
$\Symn$ acts on the irreducible representations of $H^n$ by
permuting the $n$ factors.  Let $V$ be any irreducible
representation of $H^n$, and let $G_V \subseteq \Symn$ be the
subgroup that fixes $V$ (in the action on such representations).
Then $V$ extends to a representation of $G_V \ltimes H^n$.  Taking
the tensor product with an irreducible representation $W$ of $G_V$
(with $H^n$ acting trivially on $W$) and inducing to $\Symn
\ltimes H^n$ yields an irreducible representation
$$
\Ind_{G_V \ltimes H^n}^{\Symn \ltimes H^n} (W \otimes_\C V)
$$
of $\Symn \ltimes H^n$.  All irreducible representations of $\Symn
\ltimes H^n$ arise in this way.  Two such representations are
isomorphic iff the choices of $W$ are isomorphic and the choices
of $V$ are equivalent under the action of $\Symn$ on irreducible
representations of $H^n$.  In other words, if we look at all
choices of $V$ and $W$, each representation is counted $n!/|G_V|$
times (because that is the size of $V$'s orbit under $\Symn$).

The dimension of this representation is $(n!/|G_V|) \dim(W)
\dim(V)$. Thus, because $\dim(W) \le |G_V|$ and $\sum_W \dim(W)^2
= |G_V|$,
\begin{eqnarray*}
\sum_j c_j^\omega &=& \sum_V \frac{|G_V|}{n!} \sum_W
\left(\frac{n!}{|G_V|}\dim(W)\dim(V)\right)^\omega\\
& \le & (n!)^{\omega-1}\sum_V \dim(V)^\omega = (n!)^{\omega - 1}
\left (\sum_k d_k^\omega \right )^n,
\end{eqnarray*}
as desired.
\end{proof}
}

\SubSection{Background}

In this subsection we summarize the necessary definitions and
results from \cite{CU}.

If $S$ is a subset of a group, let $Q(S)$ denote the right
quotient set of $S$, i.e., $ Q(S) = \{s_1 s_2^{-1} : s_1,s_2 \in
S\}. $

\begin{definition}[\cite{CU}]
\label{definition:realize}
A group $G$ {\em realizes} $\langle
n_1,n_2,n_3 \rangle$ if there are subsets $S_1,S_2,S_3 \subseteq
G$ such that $|S_i| = n_i$, and for $q_i \in Q(S_i)$, if
$$
q_1q_2q_3 = 1
$$
then $q_1=q_2=q_3=1$. We call this condition on $S_1,S_2,S_3$ the
{\em triple product property}.
\end{definition}

\begin{lemma}[\cite{CU}]
\label{lemma:permute} If $G$ realizes $\langle
n_1,n_2,n_3\rangle$, then it does so for every permutation of
$n_1,n_2,n_3$.
\end{lemma}

\begin{lemma}[\cite{CU}]
\label{lemma:product} If $S_1,S_2,S_3 \subseteq G$ and
$S_1',S_2',S_3' \subseteq G'$ satisfy the triple product property,
then so do the subsets $S_1 \times S_1', S_2 \times S_2', S_3
\times S_3' \subseteq G \times G'$.
\end{lemma}

\begin{theorem}[\cite{CU}]
\label{theorem:reduction} Let $R$ be any algebra over $\C$ (not
necessarily commutative). If $G$ realizes $\langle n,m,p \rangle$,
then the number of ring operations required to multiply $n \times
m$ with $m \times p$ matrices over $R$ is at most the number of
operations required to multiply two elements of $R[G]$.
\end{theorem}

One particularly useful construction from \cite{CU} involves
permutations of the points in a triangular array.  Let
$$
\Delta_n = \{(a,b,c) \in \Z^3 : \textup{$a+b+c=n-1$ and $a,b,c \ge
0$}\}.
$$
Geometrically, these triples are barycentric coordinates for a
triangular array of points with $n$ points along each side, but it
is  more convenient to manipulate them algebraically.

For $x \in \Delta_n$, we write $x = (x_1,x_2,x_3)$.  Let $H_1$,
$H_2$, and $H_3$ be the subgroups of $\Sym(\Delta_n)$ that
preserve the first, second, and third coordinates, respectively.
Specifically,
$$
H_i = \{\pi \in \Sym(\Delta_n) : \textup{$(\pi(x))_i = x_i$ for
all $x \in \Delta_n$}\}.
$$

\begin{theorem}[\cite{CU}]
\label{theorem:pseudoexponent2} The subgroups $H_1,H_2,H_3$
defined above satisfy the triple product property.
\end{theorem}

\begin{theorem}[\cite{CU}]
\label{theorem:bound} Suppose $G$ realizes $\langle n,m,p \rangle$
and the character degrees of $G$ are $\{d_i\}$. Then
$
(nmp)^{\omega/3} \le \sum_{i} d_i^{\omega}.
$
\end{theorem}

Combining Theorem~\ref{theorem:bound} with the fact that $\sum_i
d_i^2=|G|$ yields the following corollary, which is generally how
the theorem is applied:

\begin{corollary}[\cite{CU}]
\label{corollary:bound} Suppose $G$ realizes $\langle n,m,p
\rangle$ and has largest character degree $d$. Then
$(nmp)^{\omega/3} \le d^{\omega-2} |G|. $
\end{corollary}

\Section{Beating the sum of the cubes} \label{section:cubesum}

Suppose $G$ realizes $\langle n,m,p \rangle$ and has character
degrees $\{d_i\}$.  Theorem~\ref{theorem:bound} yields a
nontrivial bound on $\omega$ (by ruling out the possibility of
$\omega=3$) if and only if
$$
nmp > \sum_i d_i^3.
$$
Question~4.1 in \cite{CU} asks whether such a group exists.  In
this section we construct one, which shows that our methods do
indeed prove nontrivial bounds on $\omega$. The rest of the paper
is logically independent of this example, but it serves as
motivation for later constructions.

\remove{ We do not know of any construction that makes use of
small groups. We have used the computer program GAP \cite{GAP} to
verify by brute force search that no group of order less than
$128$ proves a nontrivial bound on $\omega$ using three subgroups
(as opposed to subsets).  Thus a construction must involve either
fairly sizable groups or subsets other than subgroups, and in fact
all of our constructions involve both.

The example in this section realizes matrix multiplication through
subsets other than subgroups. However, the subsets are close to
subgroups in the sense that they can be obtained from subgroups by
deleting a small number of elements. }

Let $H = \Cyc_n^3$, and let $G =  H^2 \rtimes \Cyc_2$, where
$\Cyc_2$ acts on $H^2$ by switching the two factors of $H$.  Let
$z$ denote the generator of $\Cyc_2$.  We write elements of $G$ in
the form $(a, b)z^i$, with $a, b \in H$ and $i \in \set{0,1}$.
Note that $z(a, b)z = (b, a)$.

Let $H_1, H_2, H_3$ be the three factors of $\Cyc_n$ in the
product $H = \Cyc_n^3$, viewed as subgroups of $H$.  For
notational convenience, let $H_4 = H_1$.  Define subsets $S_1,
S_2, S_3 \subseteq G$ by $ S_i = \set{(a,b)z^j : a \in
H_i\setminus\set{0}, b \in H_{i+1}, j \in \set{0,1}}. $ We will
prove in Lemma~\ref{lemma:simpletpp} that these subsets satisfy
the triple product property.

To analyze this construction we need very little
representation-theoretic information.  The character degrees of
$G$ are all at most $2$, because $H^2$ is an abelian subgroup of
index $2$. Then since the sum of the squares of the character
degrees is $|G|$, the sum of their cubes is at most $2|G|$, which
equals $4n^6$.

On the other hand, $|S_i| = 2n(n-1)$, so $|S_1||S_2||S_3| =
8n^3(n-1)^3$. For $n \ge 5$, this product is larger than $4n^6$.
By Corollary~\ref{corollary:bound}, $ (2n(n-1))^\omega \le
2^{\omega-2} 2n^6. $ The best bound on $\omega$ is achieved by
setting $n=17$, in which case we obtain $\omega < 2.9088$.

\remove{ It is a straightforward calculation in representation
theory to determine how many of the character degrees are $2$ (and
how many are $1$).  That can be used to improve this analysis, but
the bound on $\omega$ changes by less than $10^{-4}$, so we do not
present the details here. }

All that remains is to prove the triple product property:

\begin{lemma} \label{lemma:simpletpp}
$S_1$, $S_2$, and $S_3$ satisfy the triple product property.
\end{lemma}

\begin{proof}
Consider the triple product $q_1q_2q_3$ with $q_i \in Q(S_i)$, and
suppose it equals the identity.  Each quotient $q_i$ is either of
the form $(a_i,b_i)(-a'_i,-b'_i)$ or of the form
$(a_i,b_i)z(-a'_i,-b'_i)$, with $a_i, a'_i \in H_i$ and $b_i, b'_i
\in H_{i+1}$. There must be an even number of factors of $z$ among
the three elements $q_1, q_2, q_3$.

First, suppose there are none.  We can write $q_1q_2q_3$ as
$$
(a_1,b_1)(-a'_1, -b'_1)(a_2,b_2)(-a'_2,
-b'_2)(a_3,b_3)(-a'_3,-b'_3),
$$
where $a_i, a'_i \in H_i$ and $b_i, b'_i \in H_{i+1}$. The product
is thus equal to
$$
(a_1- a'_1 + a_2- a'_2 + a_3- a'_3, b_1- b'_1 + b_2- b'_2 + b_3-
b'_3),
$$
which is the identity iff $q_1 = q_2 = q_3 = 1$, since the triple
product property holds (trivially) for $H_1, H_2, H_3$ in $H$.

Second, suppose two of $q_1,q_2,q_3$ contain a $z$.  The product
$q_1q_2q_3$ can be simplified as above to yield a sum in each
coordinate, except now $a_i$ and $a'_i$ contribute to different
coordinates when $q_i$ contains a $z$, as do $b_i$ and $b'_i$.
There are thus two $i$'s such that $a_i$ and $a'_i$ contribute to
different coordinates.  For one of those two $i$'s, $b_{i-1}$ and
$b'_{i-1}$ contribute to the same coordinate (where we interpret
the subscripts modulo~$3$).  The sum in the other coordinate
contains one of $a_i$ and $a'_i$ but neither of $b_{i-1}$ and
$b'_{i-1}$, and thus only one summand from $H_i$ (because for each
$j$, $a_j,a'_j \in H_j$ and $b_j,b'_j \in H_{j+1}$). Since $a_i$
and $a'_i$ are nonzero by the definition of $S_i$, the product
$q_1q_2q_3$ cannot be the identity.
\end{proof}

\Section{Uniquely solvable puzzles} \label{section:usp}

In this section we define a combinatorial object called a
\textit{strong USP\/}, which gives rise to a systematic
construction of sets satisfying the triple product property in a
wreath product. Using strong USPs we achieve $\omega < 2.48$, and
we conjecture that there exist strong USPs that prove $\omega=2$.

\SubSection{USPs and strong USPs}

A \textit{uniquely solvable puzzle\/} (USP) of width $k$ is a
subset $U \subseteq \{1,2,3\}^k$ satisfying the following
property:
\begin{enumerate}
\item[] For all permutations $\pi_1,\pi_2,\pi_3 \in \Sym(U)$,
either $\pi_1=\pi_2=\pi_3$ or else there exist $u \in U$ and $i
\in [k]$ such that at least two of $(\pi_1(u))_i = 1$,
$(\pi_2(u))_i = 2$, and $(\pi_3(u))_i = 3$ hold.
\end{enumerate}

The motivation for the name ``uniquely solvable puzzle'' is that a
USP can be thought of as a jigsaw puzzle.  The puzzle pieces are
the sets $\{i : u_i=1\}$, $\{i : u_i=2\}$, and $\{i : u_i=3\}$
with $u \in U$, and the puzzle can be solved by permuting these
types of pieces according to $\pi_1$, $\pi_2$, and $\pi_3$,
respectively, and reassembling them without overlap into triples
consisting of one piece of each of the three types. The definition
requires that the puzzle must have a unique solution.

A \textit{strong USP\/} is a USP in which the defining property is
strengthened as follows:
\begin{enumerate}
\item[] For all permutations $\pi_1,\pi_2,\pi_3 \in \Sym(U)$,
either $\pi_1=\pi_2=\pi_3$ or else there exist $u \in U$ and $i
\in [k]$ such that \textit{exactly\/} two of $(\pi_1(u))_i = 1$,
$(\pi_2(u))_i = 2$, and $(\pi_3(u))_i = 3$ hold.
\end{enumerate}

One convenient way to depict USPs is by labelling a grid in which
the rows correspond to elements of the USP and the columns to
coordinates. The ordering of the rows is irrelevant.  For example,
the following labelling defines a strong USP of size $8$ and width
$6$:

\newcommand{\tcs}{\hspace{0.1cm}}
\begin{center}
\begin{tabular}{|*6{@{\tcs}c@{\tcs}|}}
\hline 3 & 3 & 3 & 3 & 3 & 3\\
\hline 1 & 3 & 3 & 2 & 3 & 3\\
\hline 3 & 1 & 3 & 3 & 2 & 3\\
\hline 1 & 1 & 3 & 2 & 2 & 3\\
\hline 3 & 3 & 1 & 3 & 3 & 2\\
\hline 1 & 3 & 1 & 2 & 3 & 2\\
\hline 3 & 1 & 1 & 3 & 2 & 2\\
\hline 1 & 1 & 1 & 2 & 2 & 2\\
\hline
\end{tabular}
\end{center}

This construction naturally generalizes as follows:

\begin{proposition} \label{proposition:easystrongusp}
For each $k \ge 1$, there exists a strong USP of size $2^k$ and
width $2k$.
\end{proposition}

\begin{proof}
Viewing $\{1,3\}^{k} \times \{2,3\}^k$ as a subset of
$\{1,2,3\}^{2k}$, we define $U$ to be
$$
\{u \in \{1,3\}^{k} \times \{2,3\}^k : \textup{for $i \in [k]$,
$u_i = 1$ iff $u_{i+k} = 2$}\}.
$$

Suppose $\pi_1,\pi_2,\pi_3 \in \Sym(U)$.  If $\pi_1 \ne \pi_3$,
then there exists $u \in U$ such that $(\pi_1(u))_i = 1$ and
$(\pi_3(u))_i = 3$ for some $i \in [k]$.  Similarly, if $\pi_2 \ne
\pi_3$, then there exists $u \in U$ such that $(\pi_2(u))_i = 2$
and $(\pi_3(u))_i = 3$ for some $i \in [2k] \setminus [k]$.  In
either case, exactly two of $(\pi_1(u))_i = 1$, $(\pi_2(u))_i =
2$, and $(\pi_3(u))_i = 3$ hold because in each coordinate only
two of the three symbols $1$, $2$, and $3$ can occur. It follows
that $U$ is a strong USP, as desired.
\end{proof}

We define the \textit{strong USP capacity\/} to be the largest
constant $C$ such that there exist strong USPs of size
$(C-o(1))^k$ and width $k$ for infinitely many values of $k$. (We
use the term ``capacity'' because this quantity is the Sperner
capacity of a certain directed hypergraph, as we explain in
Section~\ref{section:usingstpp}.) The \textit{USP capacity\/} is
defined analogously.

There is a simple upper bound for the USP capacity, which is of
course an upper bound for the strong USP capacity as well:

\begin{lemma} \label{lemma:uspcap}
The USP capacity is at most $3/2^{2/3}$.
\end{lemma}

\begin{proof}[Sketch of proof.]
A USP of width $k$ can have no repeated ``puzzle pieces'' (see the
puzzle interpretation at the beginning of this section). It
follows that the number of rows can be no larger than $O(k^2){k
\choose {k/3}} = (3/2^{2/3} + o(1))^k$.
\end{proof}

\remove{
\begin{proof}
Let $U$ be a USP of width $k$. For each triple $n_1,n_2,n_3$ of
nonnegative integers summing to $k$, define the subset
$U_{n_1,n_2,n_3}$ of $U$ to consist of all elements of $U$
containing $n_1$ entries that are $1$, $n_2$ that are $2$, and
$n_3$ that are $3$. There are $\binom{k+2}{2}$ choices of
$n_1,n_2,n_3$, so
$$
|U| \le \binom{k+2}{2} \max_{n_1,n_2,n_3} |U_{n_1,n_2,n_3}|.
$$
If two elements of $U$ have the symbol $1$ in exactly the same
locations, then letting $\pi_1$ interchange them would violate the
definition of a USP, and of course the same holds for $2$ or $3$.
Thus,
$$
|U_{n_1,n_2,n_3}| \le \min_{i} \binom{k}{n_i} \le
\left(\frac{3}{2^{2/3}}+o(1)\right)^k,
$$
where the latter inequality holds because $\min_i \binom{k}{n_i}$
is maximized when $n_1=n_2=n_3=k/3$.  It follows that $ |U| \le
\left({3}/{2^{2/3}}+o(1)\right)^k,$ as desired.
\end{proof}
}

USPs turn out to be implicit in the analysis in Coppersmith and
Winograd's paper \cite{CW}, although they are not discussed as
such. Section~6 of \cite{CW} can be interpreted as giving a
probabilistic construction showing that Lemma~\ref{lemma:uspcap}
is sharp:

\begin{theorem}[Coppersmith and Winograd \cite{CW}] \label{theorem:cw}
The USP capacity equals $3/2^{2/3}$.
\end{theorem}

We conjecture that the same is true for strong USPs:

\begin{conjecture}\label{conjecture:strongUSP}
The strong USP capacity equals $3/2^{2/3}$.
\end{conjecture}

This conjecture would imply that $\omega=2$, as we explain in the
next subsection.

\SubSection{Using strong USPs}

Given a strong USP $U$ of width $k$, let $H$ be the abelian group
of all functions from $U \times [k]$ to the cyclic group $\Cyc_m$
($H$ is a group under pointwise addition).  The symmetric group
$\Sym(U)$ acts on $H$ via
$$
\pi(h)(u,i) = h(\pi^{-1}(u),i)
$$
for $\pi \in \Sym(U)$, $h \in H$, $u \in U$, and $i \in [k]$.

Let $G$ be the semidirect product $H \rtimes \Sym(U)$, and define
subsets $S_1$, $S_2$, and $S_3$ of $G$ by letting $S_i$ consist of
all products $h \pi$ with $\pi \in \Sym(U)$ and $h \in H$
satisfying
$$
h(u,j) \ne 0 \qquad \textup{iff} \qquad u_j = i
$$
for all $u \in U$ and $j \in [k]$.

\begin{proposition} \label{proposition:goodtpp}
If $U$ is a strong USP, then $S_1$, $S_2$, and $S_3$ satisfy the
triple product property.
\end{proposition}

\begin{proof}
Consider a triple product
\begin{equation} \label{equation:triple}
h_1 \pi_1 \pi_1'^{-1} h_1'^{-1} h_2 \pi_2 \pi_2'^{-1} h_2'^{-1}
h_3 \pi_3 \pi_3'^{-1} h_3'^{-1} = 1
\end{equation}
with $h_i\pi_i, h_i'\pi_i' \in S_i$. For \eqref{equation:triple}
to hold we must have
\begin{equation} \label{equation:perm}
\pi_1 \pi_1'^{-1} \pi_2 \pi_2'^{-1} \pi_3 \pi_3'^{-1} = 1.
\end{equation}
Set $\pi = \pi_1 \pi_1'^{-1}$ and $\rho = \pi_1 \pi_1'^{-1} \pi_2
\pi_2'^{-1}$.  Then the remaining condition for
\eqref{equation:triple} to hold is that in the abelian group $H$
(with its $\Sym(U)$ action),
\begin{equation} \label{equation:tripleinH}
h_1-h_3' + \pi(h_2-h_1') + \rho(h_3-h_2') = 0.
\end{equation}
Note that
\begin{eqnarray*}
(h_1-h_3')(u,j) \ne 0 & \textup{iff} & u_j \in \{1,3\},\\
\pi(h_2-h_1')(u,j) \ne 0 & \textup{iff} &
(\pi^{-1}(u))_j \in \{2,1\}, \textup{ and}\\
\rho(h_3-h_2')(u,j) \ne 0 & \textup{iff} & (\rho^{-1}(u))_j \in
\{3,2\}.
\end{eqnarray*}
By the definition of a strong USP, either $\pi=\rho=1$ or else
there exist $u$ and $j$ such that exactly one of these three
conditions holds, in which case \eqref{equation:tripleinH} cannot
hold. Thus, $\pi = \rho = 1$, which together with
\eqref{equation:perm} implies $\pi_i=\pi'_i$ for all $i$. Then we
have $ h_1+h_2+h_3 = h_1'+h_2'+h_3', $ which implies $h_i'=h_i$
for each $i$ (because for different choices of $i$ they have
disjoint supports).  Thus, the triple product property holds.
\end{proof}

Analyzing this construction using Corollary~\ref{corollary:bound}
and the bound $[G : H] = |U|!$ on the largest character degree of
$G$ yields the following bound:

\begin{corollary} \label{corollary:omegabound}
If $U$ is a strong USP of width $k$, and $m \ge 3$ is an integer,
then $ \omega \le \frac{3\log m}{\log (m-1)} - \frac{3\log
|U|!}{|U|k \log (m-1)}. $ In particular, if the strong USP
capacity is $C$, then
$$
\omega \le \frac{3(\log m - \log C)}{\log(m-1)}.
$$
\end{corollary}

Proposition~\ref{proposition:easystrongusp} yields $\omega < 2.67$
with $m=9$. In the next subsection we prove that the strong USP
capacity is at least $2^{2/3}$ and hence $\omega < 2.48$, which is
the best bound we know how to prove using strong USPs.

If Conjecture~\ref{conjecture:strongUSP} holds, then
Corollary~\ref{corollary:omegabound} yields $\omega=2$ upon taking
$m=3$.

\SubSection{The triangle construction} \label{subsection:triangle}

The strong USP constructed in
Proposition~\ref{proposition:easystrongusp} has the property that
only two symbols (of the three possibilities $1$, $2$, and $3$)
occur in each coordinate. Every USP with this property is a strong
USP, and we can analyze exactly how large such a USP can be as
follows.

Suppose $U \subseteq \{1,2,3\}^k$ is a subset with only two
symbols occurring in each coordinate.  Let $H_1$ be the subgroup
of $\Sym(U)$ that preserves the coordinates in which only $1$ and
$2$ occur, $H_2$ the subgroup preserving the coordinates in which
only $2$ and $3$ occur, and $H_3$ the subgroup preserving the
coordinates in which only $1$ and $3$ occur.

\begin{lemma} \label{lemma:USPequiv}
The set $U$ is a USP iff $H_1$, $H_2$, and $H_3$ satisfy the
triple product property within $\Sym(U)$.
\end{lemma}

\begin{proof}
Suppose $\pi_1,\pi_2,\pi_3 \in \Sym(U)$.  The permutation $\pi_1
\pi_2^{-1}$ is not in $H_1$ iff there exists $v \in U$ and a
coordinate $i$ such that $v_i = 2$ and $((\pi_1
\pi_2^{-1})(v))_i=1$.  If we set $u = \pi_2^{-1}(v)$, then this is
equivalent to $(\pi_2(u))_i=2$ and $(\pi_1(u))_i=1$. Similarly,
$\pi_2 \pi_3^{-1} \not\in H_2$ iff there exist $u$ and $i$ such
that $(\pi_2(u))_i=2$ and $(\pi_3(u))_i=3$, and $\pi_3 \pi_1^{-1}
\not\in H_3$ iff there exist $u$ and $i$ such that
$(\pi_1(u))_i=1$ and $(\pi_3(u))_i=3$.

Thus, $U$ is a USP iff for all $\pi_1,\pi_2,\pi_3$, if $\pi_1
\pi_2^{-1} \in H_1$, $\pi_2 \pi_3^{-1} \in H_2$, and $\pi_3
\pi_1^{-1} \in H_3$, then $\pi_1=\pi_2=\pi_3$.  That is equivalent
to the triple product property for $H_1$, $H_2$, and $H_3$: recall
that because these are subgroups, the triple product property says
that for $h_i \in H_i$, $h_1h_2h_3=1$ iff $h_1=h_2=h_3=1$.  Any
three elements $h_1,h_2,h_3$ satisfying $h_1h_2h_3=1$ can be
written in the form $h_1 = \pi_1 \pi_2^{-1}$, $h_2 = \pi_2
\pi_3^{-1}$, and $h_3 = \pi_3 \pi_1^{-1}$.
\end{proof}

\begin{proposition} \label{proposition:bigstrong}
For each $k \ge 1$, there exists a strong USP of size
$2^{k-1}(2^k+1)$ and width $3k$.
\end{proposition}

It follows that the strong USP capacity is at least $2^{2/3}$ and
$\omega < 2.48$.

\begin{proof}
Consider the triangle
$$
\Delta_{n} = \{(a,b,c) \in \Z^3 : \textup{$a+b+c=n-1$ and $a,b,c
\ge 0$}\},
$$
with $n = 2^k$, and let $H_1$, $H_2$, and $H_3$ be the subgroups
of $\Sym(\Delta_n)$ preserving the first, second, and third
coordinates, respectively.  By
Theorem~\ref{theorem:pseudoexponent2}, these subgroups satisfy the
triple product property in $\Sym(\Delta_n)$.

To construct the desired strong USP, choose a subset $U \subseteq
\{1,2,3\}^{3k}$ as follows.  Among the first $k$ coordinates, only
$1$ and $2$ will occur, among the second $k$ only $2$ and $3$, and
among the third $k$ only $1$ and $3$.  In each of these three
blocks of $k$ coordinates, there are $2^k$ possible patterns that
be made using the two available symbols.  Number these patterns
arbitrarily from $0$ to $2^k-1$ (each number will be used for
three patterns, one for each pair of symbols).  The elements of
$U$ will correspond to elements of $\Delta_n$. In particular, the
element of $U$ corresponding to $(a,b,c) \in \Delta_n$ will have
the $a$-th pattern in the first $k$ coordinates, the $b$-th in the
second $k$, and the $c$-th in the third.  It follows from
Lemma~\ref{lemma:USPequiv} that $U$ is a strong USP.
\end{proof}

One can show using Lemma~\ref{lemma:USPequiv} that this
construction is optimal:

\begin{corollary}
If $U$ is a USP of width $k$ such that only two symbols occur in
each coordinate, then $ |U| \le (2^{2/3}+o(1))^k. $
\end{corollary}

The condition of using only two symbols in each coordinate is
highly restrictive, but we have been unable to improve on
Proposition~\ref{proposition:bigstrong}.  However, we know of no
upper bound on the size of a strong USP besides
Lemma~\ref{lemma:uspcap}, and we see no reason why
Conjecture~\ref{conjecture:strongUSP} should not be true.

\Section{The simultaneous double product property}
\label{section:sdpp}

There are at least two natural avenues for improving the
construction from Subsection~\ref{subsection:triangle}.  In the
combinatorial direction, one might hope to replace the strong USP
of Proposition~\ref{proposition:bigstrong} with a larger one; this
will reach exponent~$2$ if Conjecture~\ref{conjecture:strongUSP}
holds. In the algebraic direction, one might hope to keep the
combinatorial structure of the triangle construction in place
while modifying the underlying group. Such a modification can be
carried out using the \textit{simultaneous double product
property\/} defined below, and we conjecture that it reaches
$\omega=2$ as well (Conjecture~\ref{conjecture:sdpp}).

We say that subsets $S_1, S_2$ of a group $H$ satisfy the {\em
double product property} if
\[q_1q_2 = 1 \qquad \textup{implies} \qquad q_1 = q_2 = 1,\]
where $q_i \in Q(S_i)$.

\begin{definition} \label{definition:sdpp}
We say that $n$ pairs of subsets $A_i, B_i$ (for $1 \le i \le n$)
of a group $H$ satisfy the {\em simultaneous double product
property} if
\begin{itemize}
\item for all $i$, the pair $A_i, B_i$ satisfies the double product property,
and

\item for all $i, j, k$,
\[a_i(a_j')^{-1}b_j(b_k')^{-1} = 1 \qquad \textup{implies}\qquad i = k,\]
where $a_i \in A_i$, $a_j' \in A_j$, $b_j \in B_j$, and $b_k' \in
B_k$.
\end{itemize}
\end{definition}

A convenient reformulation is that if one looks at the sets
$$
A_i^{-1}B_j = \{a^{-1}b : a \in A_i, b \in B_j\},
$$
those with $i=j$ are disjoint from those with $i \ne j$.

For a trivial example, set $H = \Cyc_n^k \times \Cyc_n$, and set
$A_i = \{(x,i): x \in \Cyc_n^k\}$ and $B_i = \{(0,i)\}$. Then the
pairs $A_i,B_i$ for $i \in \Cyc_n$ satisfy the simultaneous double
product property.

\begin{lemma} \label{lemma:sdppprod}
If $n$ pairs of subsets $A_i, B_i \subseteq H$ satisfy the
simultaneous double product property, and $n'$ pairs of subsets
$A_i', B_i' \subseteq H'$ satisfy the simultaneous double product
property, then so do the $nn'$ pairs of subsets $A_i \times A_j',
B_i \times B_j' \subseteq H \times H'$.
\end{lemma}

Pairs $A_i, B_i$ satisfying the simultaneous double product
property in group $H$ can be transformed into subsets satisfying
the triple product property via a construction similar to the one
in Section~\ref{section:usp}. Recall that
$$
\Delta_n = \{(a,b,c) \in \Z^3 : \textup{$a+b+c=n-1$ and $a,b,c \ge
0$}\}.
$$
Given $n$ pairs of subsets $A_i, B_i$ in $H$ for $0 \le i \le
n-1$, we define triples of subsets in $H^3$ indexed by $v =
(v_1,v_2,v_3) \in \Delta_n$ as follows:
\begin{eqnarray*}
\widehat{A}_v & = & A_{v_1} \times \{1\} \times B_{v_3} \\
\widehat{B}_v & = & B_{v_1} \times A_{v_2} \times \{1\} \\
\widehat{C}_v & = & \{1\} \times B_{v_2} \times A_{v_3}
\end{eqnarray*}

\begin{theorem}
\label{theorem:SDPP2TPP} If $n$ pairs of subsets $A_i, B_i
\subseteq H$ (with $0 \le i \le n-1$) satisfy the simultaneous
double product property, then the following subsets $S_1, S_2,
S_3$ of $G = (H^3)^{\Delta_n} \rtimes \Sym(\Delta_n)$ satisfy the
triple product property:
\begin{eqnarray*}
S_1 & = & \{\widehat{a}\pi: \pi \in \Sym(\Delta_n), \widehat{a}_v
\in \widehat{A}_v \mbox{ for all $v$}\} \\
S_2 & = & \{\widehat{b}\pi: \pi \in \Sym(\Delta_n), \widehat{b}_v
\in \widehat{B}_v \mbox{ for all $v$}\} \\
S_3 & = & \{\widehat{c}\pi: \pi \in \Sym(\Delta_n), \widehat{c}_v
\in \widehat{C}_v \mbox{ for all $v$}\}
\end{eqnarray*}
\end{theorem}

The proof uses Theorem~\ref{theorem:pseudoexponent2} and is
similar to the proof of Proposition~\ref{proposition:goodtpp}; it
can be found in the full version of this paper.

\begin{theorem}
\label{theorem:SDPP-bound} If $H$ is a finite group with character
degrees $\{d_k\}$, and $n$ pairs of subsets $A_i, B_i \subseteq H$
satisfy the simultaneous double product property, then
\[\sum_{i=1}^n(|A_i||B_i|)^{\omega/2} \le \left (\sum_k d_k^\omega
\right)^{3/2}.\]
\end{theorem}

Using this theorem, the example after
Definition~\ref{definition:sdpp} recovers the trivial bound
$\omega \le 3$ as $k \to \infty$.

\begin{proof}[Proof of Theorem \ref{theorem:SDPP-bound}]
Let $A_i', B_i'$ be the $N$-fold direct product of the pairs $A_i,
B_i$ via Lemma~\ref{lemma:sdppprod}, and let $\mu$ be an arbitrary
$n$-vector of nonnegative integers for which $\sum_{i=1}^n \mu_i =
N$. Among the pairs $A_i', B_i'$ are $M = \binom{N}{\mu}$ pairs
for which $|A_i'||B_i'| = \prod_{i = 1}^n (|A_i||B_i|)^{\mu_i};$
call this quantity $L$. Set $P = |\Delta_M|$, so $P = M(M+1)/2$.
The three subsets in Theorem~\ref{theorem:SDPP2TPP} each have size
$P!L^{P}$. By Theorem~\ref{theorem:bound} and Lemma
\ref{lemma:wreath-char-degrees} we obtain $(P!L^P)^\omega \le
(P!)^{\omega - 1}\left ( \sum_k d_k^\omega \right )^{3NP}.$ Taking
$2P$-th roots and letting $N \rightarrow \infty$ yields
$$
\binom{N}{\mu}\left (\prod_{i = 1}^n (|A_i||B_i|)^{\mu_i}\right
)^{\omega/2} \le \left (\sum_k d_k^\omega \right )^{3N/2}.
$$
Finally, we apply Lemma \ref{lemma:geom->arith} with $s_i
=(|A_i||B_i|)^{\omega/2}$ and $C = (\sum_k d_k^\omega)^{3/2}$ to
obtain the stated inequality.
\end{proof}

It is convenient to use two parameters $\alpha$ and $\beta$ to
describe pairs satisfying the simultaneous double product
property: if there are $n$ pairs, choose $\alpha$ and $\beta$ so
that $|A_i||B_i| \ge n^\alpha$ for all $i$ and $|H|=n^\beta$. If
$H$ is abelian Theorem~\ref{theorem:SDPP-bound} implies $\omega
\le (3\beta-2)/\alpha$.

The best construction we know is the following:

\begin{proposition}
For each $m \ge 2$, there is a construction in $\Cyc_m^{2\ell}$
satisfying the simultaneous double product property with $\alpha =
\log_2 (m-1) +o(1)$ and $\beta = \log_2 m +o(1)$ as $\ell \to
\infty$.
\end{proposition}

Taking $m=6$ yields exactly the same bound as in
Subsection~\ref{subsection:triangle} ($\omega < 2.48$).

\begin{proof}
Let $n = \binom{2\ell}{\ell}$.  Then $n=2^{2\ell(1-o(1))}$ so
$\beta = \log_2 m + o(1)$. For each subset $S$ of the $2\ell$
coordinates of $\Cyc_m^{2\ell}$ with $|S|=\ell$, let $A_S$ be the
set of elements that are nonzero in those coordinates and zero in
the others.  Let $\overline{S}$ denote the complement of $S$, and
set $B_S = A_{\overline{S}}$.  For each $S$, we have $|A_S||B_S| =
(m-1)^{2\ell}$, so $\alpha = \log_2(m-1) + o(1)$.

We will show that the pairs $A_S,B_S$ satisfy the simultaneous
double product property. Each pair $A_S, B_S$ clearly satisfies
the double product property, because the elements of $A_S$ and
$B_S$ are supported on disjoint sets of coordinates. Each element
of $B_S - A_S$ is nonzero in every coordinate, but if $Q \ne R$
then there is a coordinate in $\overline{R}\cap Q$ (note that this
is why we require $|Q|=|R|$). Each element of $B_Q -A_R$ vanishes
in that coordinate, so
$$
(B_Q - A_R) \cap (B_S - A_S) =\emptyset
$$
as desired.
\end{proof}

The only limitations we know of on the possible values of $\alpha$
and $\beta$ are the following:

\begin{proposition} \label{proposition:alphabeta}
If $n$ pairs of subsets $A_i, B_i \subseteq H$ satisfy the
simultaneous double product property, with $|A_i||B_i| \ge
n^\alpha$ for all $i$ and $|H|=n^\beta$, then $\alpha \le \beta$
and $\alpha+2 \le 2\beta$.
\end{proposition}

\begin{proof}[Sketch of proof.]
The inequality $\alpha \le \beta$ follows immediately from the
double product property; the second inequality is proved using the
fact that the sets $A_1,\dots,A_n$ are pairwise disjoint, as are
$B_1,\dots,B_n$.
\end{proof}

\remove{
\begin{proof}
The inequality $\alpha \le \beta$ follows immediately from the
double product property, since that means that the quotient map
$(a,b) \mapsto a^{-1}b$ from $A_i \times B_i$ to $H$ is injective.

For the other inequality, first note that $A_1,\dots,A_n$ are
disjoint (if $x \in A_i \cap A_j$ with $i \ne j$, and $y \in B_i$,
then $x^{-1}y \in (A_i^{-1}B_i)\cap(A_j^{-1}B_i)$, which is
impossible). Similarly, $B_1,\dots,B_n$ are also disjoint.  It
follows that the map
$$
\Symn \times \Symn \times \prod_{i=1}^n A_i \times \prod_{i=1}^n
B_i \to (H^n)^2
$$
defined by $(\pi,\rho,a,b) \mapsto (\pi a,\rho b)$ is injective.
Here, the group $\Symn$ acts by permuting the $n$ coordinates.
Comparing the sizes of these sets yields $ (n!)^2 (n^\alpha)^n \le
(n^\beta)^{2n},$ which implies $2\beta \ge \alpha+2$ as $n
\to\infty$. Note that by taking direct powers via
Lemma~\ref{lemma:sdppprod}, one can take $n$ arbitrarily large
without changing $\alpha$ and $\beta$.
\end{proof}
}

The most important case is when $H$ is an abelian group.  There
the bound on $\omega$ is $\omega \le (3\beta-2)/\alpha$, and
Proposition~\ref{proposition:alphabeta} shows that the only way to
achieve $\omega=2$ is $\alpha=\beta=2$.  We conjecture that that
is possible:

\begin{conjecture} \label{conjecture:sdpp}
For arbitrarily large $n$, there exists an abelian group $H$ with
$n$ pairs of subsets $A_i,B_i$ satisfying the simultaneous double
product property such that $|H| = n^{2+o(1)}$ and $|A_i||B_i| \ge
n^{2-o(1)}$.
\end{conjecture}

\Section{The simultaneous triple product property}
\label{section:stpp}

Each of our constructions of a group proving a nontrivial bound on
$\omega$ has the same general form, namely a semidirect product of
a permutation group with an abelian group.  The crucial part of
such a construction is the way in which the abelian part is
apportioned among the three subsets satisfying the triple product
property.

This apportionment can be viewed as reducing several independent
matrix multiplication problems to a single group algebra
multiplication, using triples of subsets satisfying the
simultaneous triple product property:

\begin{definition}
We say that $n$ triples of subsets $A_i, B_i, C_i$ (for $1 \le i
\le n$) of a group $H$ satisfy the {\em simultaneous triple
product property} if
\begin{itemize}
\item for each $i$, the three subsets $A_i, B_i, C_i$ satisfy the
triple product property, and

\item for all $i, j, k$,
$$
a_i (a_j')^{-1} b_j (b_k')^{-1} c_k (c_i')^{-1} = 1 \quad
\textup{implies} \quad i = j = k
$$
for $a_i \in A_i$, $a_j' \in A_j$, $b_j \in B_j$, $b_k' \in B_k$,
$c_k \in C_k$ and $c_i' \in C_i$.
\end{itemize}
We say that such a group {\em simultaneously realizes} $\langle
|A_1|,|B_1|,|C_1| \rangle, \dots, \langle
|A_n|,|B_n|,|C_n|\rangle$.
\end{definition}

In most applications the group $H$ will be abelian, in which case
it is more conventional to use additive notation.  In this
notation the implication above becomes
$$
a_i - a_j' + b_j - b_k' + c_k - c_i' = 0 \qquad \textup{implies}
\qquad i = j = k.
$$

As an example, let $H = \Cyc_n^3$, and call the three factors
$H_1$, $H_2$, and $H_3$.  Define the following sets:
\begin{eqnarray*}
A_1 = H_1 \setminus \{0\}, \quad B_1 = H_2 \setminus \{0\}, \quad
C_1 = H_3 \setminus \{0\} \\
A_2 = H_2 \setminus \{0\}, \quad B_2 = H_3 \setminus \{0\}, \quad
C_2 = H_1 \setminus \{0\}
\end{eqnarray*}
This construction is based on the one in
Section~\ref{section:cubesum}, except that this one is slightly
more symmetrical.

\begin{proposition} \label{proposition:stppexample}
The two triples $A_1,B_1,C_1$ and $A_2,B_2,C_2$ satisfy the
simultaneous triple product property.
\end{proposition}

\begin{proof}[Sketch of proof.]
Each triple clearly satisfies the triple product property in
isolation, so we need only deal with the second condition in the
definition. For $i \in \{1,2\}$ define $U_i = A_i - C_i$, $V_i =
B_i-A_i$, and $W_i = C_i-B_i$.  The theorem follows from arguing
that if $u_i+v_j+w_k=0$ (with $u_i \in U_i$, $v_j \in V_j$, and
$w_k \in W_k$) then $i=j=k$.
\end{proof}

\remove{
\begin{proof}
Each triple clearly satisfies the triple product property in
isolation, so we need only deal with the second condition in the
definition. For $i \in \{1,2\}$ define $U_i = A_i - C_i$, $V_i =
B_i-A_i$, and $W_i = C_i-B_i$.  What we must prove is that if
$u_i+v_j+w_k=0$ with $u_i \in U_i$, $v_j \in V_j$, and $w_k \in
W_k$, then $i=j=k$.

We have
\begin{eqnarray*}
U_1 = W_2 &=& \{(x,0,z) \in \Cyc_n^3 : x \ne 0, z \ne 0\}, \\
V_1 = U_2 &=& \{(x,y,0) \in \Cyc_n^3 : x \ne 0, y \ne 0\}, \\
W_1 = V_2 &=& \{(0,y,z) \in \Cyc_n^3 : y \ne 0, z \ne 0\}.
\end{eqnarray*}
If $i$, $j$, and $k$ are not all equal, then two of them must be
equal but different from the third.  In each such case, $U_i$,
$V_j$, and $W_k$ comprise exactly two of the three subsets of
$\Cyc_n^3$ defined in the equations above, with one of those two
sets occurring twice.  The sum of an element from each cannot
vanish, since in the repeated set one coordinate is zero, and the
other set is always nonzero in that coordinate.
\end{proof}
}

The reason for the strange condition in the definition of the
simultaneous triple product property is that it is exactly what is
needed to reduce several independent matrix multiplications to one
group algebra multiplication.

\begin{theorem}
\label{theorem:sttpreduction} Let $R$ be any algebra over $\C$. If
$H$ simultaneously realizes $\langle n_1,m_1,p_1
\rangle,\dots,\langle n_k,m_k,p_k \rangle$, then the number of
ring operations required to perform $k$ independent matrix
multiplications of sizes $n_1 \times m_1$ by $m_1 \times p_1$,
\dots, $n_k \times m_k$ by $m_k \times p_k$ is at most the number
of operations required to multiply two elements of $R[H]$.
\end{theorem}

The proof is similar to that of Theorem~\ref{theorem:reduction}:

\begin{proof}
Suppose $H$ simultaneously realizes $\langle n_1,m_1,p_1 \rangle$,
\dots, $\langle n_k,m_k,p_k \rangle$ via triples $N_i,M_i,P_i$
with $1 \le i \le k$.  Let $A_i$ be an $n_i \times m_i$ matrix and
$B_i$ an $m_i \times p_i$ matrix.  We will index the rows and
columns of $A_i$ with the sets $N_i$ and $M_i$, respectively,
those of $B_i$ with $M_i$ and $P_i$, and those of $A_iB_i$ with
$N_i$ and $P_i$.

Consider the product of these two elements of $R[H]$:
\begin{eqnarray*}
\sum_{i=1}^k\sum_{s\in N_i, t \in M_i} (A_i)_{s,t} s^{-1} t; & &
\sum_{i=1}^k\sum_{t' \in M_i, u \in P_i} (B_i)_{t',u} t'^{-1} u.
\end{eqnarray*}

We have $ (s^{-1} t) (t'^{-1} u) = s'^{-1} u' $ with $s \in N_i$,
$t \in M_i$, $t' \in M_j$, $u \in P_j$, $s' \in N_k$, and $u' \in
P_k$ iff $i=j=k$ and $s=s'$, $t=t'$, and $u=u'$, so the
coefficient of $s^{-1} u$ in the product is $ \sum_{t \in T}
(A_i)_{s,t} (B_i)_{t,u} = (A_iB_i)_{s,u}. $ Thus, one can simply
read off the matrix products from the group algebra product by
looking at the coefficients of $s^{-1} u$ with $s \in N_i, u \in
P_i$, and the theorem follows.
\end{proof}

Other results about the triple product property also generalize
straightforwardly to the simultaneous triple product property,
such as the following lemma:

\begin{lemma} \label{lemma:directprod}
If $n$ triples of subsets $A_i, B_i, C_i \subseteq H$ satisfy the
simultaneous triple product property, and $n'$ triples of subsets
$A_j', B_j', C_j' \subseteq H'$ satisfy the simultaneous triple
product property, then so do the $nn'$ triples of subsets $A_i
\times A_j', B_i \times B_j', C_i \times C_j' \subseteq H \times
H'$.
\end{lemma}

By Sch\"onhage's asymptotic sum inequality ((15.11) in
\cite{BCS}), one can deduce a bound on $\omega$ from the
simultaneous triple product property:

\begin{theorem} \label{theorem:asi}
If a group $H$ simultaneously realizes $\langle
a_1,b_1,c_1\rangle,\dots,\langle a_n,b_n,c_n \rangle$ and has
character degrees $\{d_k\}$, then
$
\sum_{i=1}^n (a_ib_ic_i)^{\omega/3} \le \sum_k d_k^\omega.
$
\end{theorem}

Frequently $H$ will be abelian, in which case $\sum_k d_k^\omega =
|H|$. That occurs in the example from
Proposition~\ref{proposition:stppexample}, which proves that
$\omega < 2.93$ using Theorem~\ref{theorem:asi}.

In Section~\ref{section:wreath} we provide a proof of
Theorem~\ref{theorem:asi} completely within our group-theoretic
framework, and show furthermore that any bound on $\omega$ that
can be achieved using the simultaneous triple product property can
also be achieved using the ordinary triple product property. Thus,
there is no added generality from the simultaneous triple product
property, but it is an important organizing principle.

\Section{Using the simultaneous triple product property}
\label{section:usingstpp}

Every construction we have found of a group proving a nontrivial
bound on $\omega$ has at its core a simultaneous triple product
property construction in an abelian group.  Each construction also
involves a wreath product, but as explained in
Section~\ref{section:wreath} that is a general tool for dealing
with the simultaneous triple product property. Given
Theorem~\ref{theorem:asi}, which can be proved either via the
wreath product construction of Section~\ref{section:wreath} or
using the asymptotic sum inequality, one can dispense with
non-abelian groups entirely. In this section we explain how to
interpret each of our constructions in this setting.

\SubSection{Local strong USPs}
\label{subsection:loc-strong-usp}

A \textit{local strong USP} of width $k$ is a subset $U \subseteq
\{1,2,3\}^k$ such that for each ordered triple $(u,v,w) \in U^3$,
with $u$, $v$, and $w$ not all equal, there exists $i \in [k]$
such that $(u_i,v_i,w_i)$ is an element of
$$
\{(1,2,1),(1,2,2),(1,1,3),(1,3,3),(2,2,3),(3,2,3)\}.
$$

\begin{lemma} \label{lemma:localisglobal}
Every local strong USP is a strong USP.
\end{lemma}

\remove{
\begin{proof}
Let $U$ be a local strong USP, and suppose $\pi_1,\pi_2,\pi_3 \in
\Sym(U)$. If $\pi_1$, $\pi_2$, and $\pi_3$ are not all equal, then
there exists $u \in U$ such that $\pi_1(u)$, $\pi_2(u)$, and
$\pi_3(u)$ are not all equal.  There exists $i \in [k]$ such that
$((\pi_1(u))_i,(\pi_2(u))_i,(\pi_3(u))_i)$ is in $
\{(1,2,1),(1,2,2),(1,1,3),(1,3,3),(2,2,3),(3,2,3)\}, $ and hence
exactly two of $(\pi_1(u))_i = 1$, $(\pi_2(u))_i = 2$, and
$(\pi_3(u))_i = 3$ hold, as desired.
\end{proof}
}

The proof is straightforward, and omitted from this version of the
paper.

The reason for the word ``local'' is that local strong USPs
satisfy a condition for every triple of rows, rather than a weaker
global condition on permutations. The advantage of local strong
USPs is that they lead naturally to a construction satisfying the
simultaneous triple product property:

\begin{theorem} \label{theorem:SUSPtoSTPP}
Let $U$ be a local strong USP of width $k$, and for each $u \in U$
define subsets $A_u,B_u,C_u \subseteq \Cyc_\ell^k$ by
\begin{eqnarray*}
A_u &=& \{x \in \Cyc_\ell^{k} : x_j \ne 0 \textup{ iff $u_j = 1$}\},\\
B_u &=& \{x \in \Cyc_\ell^{k} : x_j \ne 0 \textup{ iff $u_j = 2$}\}, \textup{ and}\\
C_u &=& \{x \in \Cyc_\ell^{k} : x_j \ne 0 \textup{ iff $u_j =
3$}\}.
\end{eqnarray*}
Then the triples $A_u,B_u,C_u$ satisfy the simultaneous triple
product property.
\end{theorem}

Note that this construction isolates the key idea behind
Proposition~\ref{proposition:goodtpp}.

\begin{proof}
Suppose $u,v,w \in U$ are not all equal and
$$
a_u  - a_v' + b_v  - b_w' + c_w - c_u'  =  0
$$
with $a_u \in A_u$, $a_v' \in A_v$, $b_v \in B_v$, $b_w' \in B_w$,
$c_w \in C_w$ and $c_u' \in C_u$. By the definition of a local
strong USP, there exists $i \in [k]$ such that $(u_i,v_i,w_i)$ is
in
$$
\{(1,2,1),(1,2,2),(1,1,3),(1,3,3),(2,2,3),(3,2,3)\}.
$$
In each of these cases exactly one of $a_u,a_v',b_v,b_w',c_w,c_u'$
is nonzero, namely $a_v'$, $b_v$, $a_u$, $c_u'$, $b_w'$, and
$c_w$, respectively.  Thus, in each case the equation $ a_u  + b_v
+ c_w  =  a_v'+b_w'+c_u' $ is impossible, so $u=v=w$, as desired.

All that remains is to show that for each $u$, the sets
$A_u,B_u,C_u$ satisfy the triple product property, which is
trivial (they are supported on disjoint sets of coordinates).
\end{proof}

At first glance the definition of a local strong USP appears far
stronger than that of a strong USP.  For example, the strong USPs
constructed in Subsection~\ref{subsection:triangle} are not local
strong USPs.  However, it turns out that any bound on $\omega$
that can be proved using strong USPs can be proved using local
strong USPs:

\begin{proposition} \label{proposition:globaltolocal}
The strong USP capacity is achieved by local strong USPs.  In
particular, given any strong USP $U$ of width $k$, there exists a
local strong USP of size $|U|!$ and width $|U|k$.
\end{proposition}

\begin{proof}
Let $U$ be a strong USP of width $k$, and fix an arbitrary
ordering $u_1,u_2,\dots,u_{|U|}$ of the elements of $U$. For each
$\pi \in \Sym(U)$, let $U_\pi \in \{1,2,3\}^{|U|k}$ be the
concatenation of $\pi(u_1),\pi(u_2),\dots,\pi(u_{|U|})$.  Then the
set of all vectors $U_\pi$ is a local strong USP: given any three
elements $U_{\pi_1}$, $U_{\pi_2}$, and $U_{\pi_3}$ with
$\pi_1,\pi_2,\pi_3$ not all equal, by the definition of a strong
USP there exist $u \in U$ and $i \in [k]$ such that exactly two of
$(\pi_1(u))_i = 1$, $(\pi_2(u))_i = 2$, and $(\pi_3(u))_i = 3$
hold.  Then in the coordinate indexed by $u$ and $i$, the vectors
$U_{\pi_1}$, $U_{\pi_2}$, and $U_{\pi_3}$ have entries among
$(1,2,1)$, $(1,2,2)$, $(1,1,3)$, $(1,3,3)$, $(2,2,3)$, $(3,2,3)$,
as desired.
\end{proof}

Proposition~\ref{proposition:globaltolocal} explains the choice of
the word ``capacity'': optimizing the size of a local strong USP
amounts to determining the Sperner capacity of a certain directed
hypergraph (see \cite{Simonyi} for background on Sperner
capacity).  The full version of this paper will explain this
perspective more completely.

\SubSection{Triangle-free sets}

The construction in Theorem~\ref{theorem:SDPP2TPP} is also easily
interpreted in terms of the simultaneous triple product property.
Recall the construction of triples $\widehat{A}_v, \widehat{B}_v,
\widehat{C}_v$ indexed by $v \in \Delta_n$, defined before
Theorem~\ref{theorem:SDPP2TPP}. These triples almost satisfy the
simultaneous triple product property, in the following sense: if $
a_u (a_v')^{-1} b_v (b_w')^{-1} c_w (c_u')^{-1} = 1 $ then it
follows from the simultaneous double product property that
$u_1=w_1$, $v_2 = u_2$, and $w_3=v_3$. Call a subset $S$ of
$\Delta_n$ \textit{triangle-free\/} if for all $u,v,w \in S$
satisfying $u_1=w_1$, $v_2 = u_2$, and $w_3=v_3$, it follows that
$u=v=w$. Thus, the triples $\widehat{A}_v, \widehat{B}_v,
\widehat{C}_v$ with $v$ in a triangle-free subset of $\Delta_n$
satisfy the triple product property.

The critical question is whether there is a triangle-free subset
of $\Delta_n$ of size $|\Delta_n|^{1-o(1)}$. We give a simple
construction achieving this using Salem-Spencer sets (see
\cite{SS}). Let $T$ be a subset of $[\lfloor n/2 \rfloor]$ of size
$n^{1-o(1)}$ that contains no three-term arithmetic progression.
The following lemma is easily proved:

\begin{lemma}
The subset $\{(a,b,c) \in \Delta_n : b-a \in T\}$ is triangle-free
and has size $|\Delta_n|^{1-o(1)}$.
\end{lemma}

\SubSection{Local USPs and generalizations}

USPs also have a local version, just as strong USPs do.  A
\textit{local USP\/} is defined analogously to a local strong USP,
except that the triple $(1,2,3)$ is allowed in addition to
$(1,2,1)$, $(1,2,2)$, $(1,1,3)$, $(1,3,3)$, $(2,2,3)$, and
$(3,2,3)$.  Local USPs are USPs, and they achieve the USP
capacity; the proofs are analogous to those for
Lemma~\ref{lemma:localisglobal} and
Proposition~\ref{proposition:globaltolocal}.  In what follows we
place this construction in a far broader context:

\begin{definition}
Let $H$ be a finite abelian group.  An \emph{$H$-chart}
$\Chart=(\Gamma,A,B,C)$ consists of a finite set of symbols
$\Gamma$, together with three mappings $A,B,C \,:\, \Gamma
\rightarrow 2^H$ such that for each $x \in \Gamma,$ the sets
$A(x),B(x),C(x)$ satisfy the triple product property.  Let
$\HypGph(\Chart) \subseteq \Gamma^3$ denote the set of ordered
triples $(x,y,z)$ such that
$$
0 \not\in A(x)-A(y)+B(y)-B(z)+C(z)-C(x).
$$
A \emph{local $\Chart$-USP of width $k$} is a subset $U \subseteq
\Gamma^k$ such that for each ordered triple $(u,v,w) \in U^3$,
with $u,v,w$ not all equal, there exists $i \in [k]$ such that
$(u_i,v_i,w_i) \in \HypGph(\Chart).$
\end{definition}

For example, a local USP is a $\Chart$-USP for the
$\Cyc_\ell$-chart $\Chart = (\{1,2,3\},A,B,C)$ with $A,B,C$
defined as follows (below, $\widehat{H} = \Cyc_\ell \setminus
\{0,1\}$):
$$
\begin{array}{lll}
A(1) = \{0\} & B(1) = -\widehat{H} & C(1) = \{0\} \\
A(2) = \{1\} & B(2) = \{0\} & C(2) = \widehat{H} \\
A(3) = \widehat{H} & B(3) = \{0\} & C(3) = \{0\}
\end{array}
$$

\begin{theorem} \label{theorem:chartstpp}
Let $H$ be a finite abelian group, $\Chart$ an $H$-chart, and $U$
a local $\Chart$-USP of width $k$.  For each $u \in U$ define
subsets $A_u, B_u, C_u \subseteq H^k$ by
$$
A_u = \prod_{i=1}^k A(u_i), \quad B_u = \prod_{i=1}^k B(u_i),
\quad C_u = \prod_{i=1}^k C(u_i).
$$
Then these triples of subsets satisfy the simultaneous
triple product property.
\end{theorem}

Together with the example above, this theorem gives an analogue of
Theorem~\ref{theorem:SUSPtoSTPP} for local USPs.  Using
Theorem~\ref{theorem:cw}, this example achieves $\omega < 2.41$.

Using a more complicated chart with 24 symbols, the bound $\omega
< 2.376$ from~\cite{CW} may be derived from
Theorem~\ref{theorem:chartstpp}. For details, see the full version
of this paper.

\Section{The wreath product construction} \label{section:wreath}

It remains to prove Theorem~\ref{theorem:asi} using purely
group-theoretic means.  Besides giving a self-contained proof,
this will also show that the ordinary triple product property from
Definition~\ref{definition:realize} is as strong as the
simultaneous triple product property, in the sense that any bound
that can be derived from Theorem~\ref{theorem:asi} can be proved
using Theorem~\ref{theorem:bound} as well.

To prove Theorem~\ref{theorem:asi}, we make use of a wreath
product construction. Let $H$ be a group, and define $G = \Symn
\ltimes H^n$, where the symmetric group $\Symn$ acts on $H^n$ from
the right by permuting the coordinates according to $(h^\pi)_i =
h_{\pi(i)}$. We write elements of $G$ as $h \pi$ with $h \in H^n$
and $\pi \in \Symn$.

\begin{theorem}
\label{theorem:STPP2TPP} If $n$ triples of subsets $A_i, B_i, C_i
\subseteq H$ satisfy the simultaneous triple product property,
then the following subsets $H_1, H_2, H_3$ of $G = \Symn \ltimes
H^n$ satisfy the triple product property:
\begin{eqnarray*}
H_1 & = & \{h \pi: \pi \in \Symn, h_i \in A_i \textup{ for each $i$}\} \\
H_2 & = & \{h \pi: \pi \in \Symn, h_i \in B_i \textup{ for each $i$}\} \\
H_3 & = & \{h \pi: \pi \in \Symn, h_i \in C_i \textup{ for each
$i$}\}
\end{eqnarray*}
\end{theorem}

\begin{proof}
The proof is analogous to that of
Proposition~\ref{proposition:goodtpp}.  Consider a triple product
\begin{equation} \label{equation:stpptriple}
h_1 \pi_1 \pi_1'^{-1} h_1'^{-1} h_2 \pi_2 \pi_2'^{-1} h_2'^{-1}
h_3 \pi_3 \pi_3'^{-1} h_3'^{-1} = 1
\end{equation}
with $h_i\pi_i, h_i'\pi_i' \in H_i$.  (Note that these subscripts
index $h_1,h_2,h_3$, rather than describing coordinates of a
single $h \in H$.  Once understood that should not cause
confusion.) For \eqref{equation:stpptriple} to hold we must have
\begin{equation} \label{equation:stppperm}
\pi_1 \pi_1'^{-1} \pi_2 \pi_2'^{-1} \pi_3 \pi_3'^{-1} = 1.
\end{equation}
Set $\pi = \pi_1 \pi_1'^{-1}$ and $\rho = \pi_1 \pi_1'^{-1} \pi_2
\pi_2'^{-1}$.  Then the remaining condition for
\eqref{equation:stpptriple} to hold is that in the group $H^n$
with its right $\Symn$ action,
$$
h_3'^{-1} h_1 \big(h_1'^{-1} h_2\big)^\pi \big(h_2'^{-1}
h_3\big)^\rho = 1
$$
In other words, for each coordinate $i$,
$$
\big(h_3'^{-1}\big)_i \big(h_1\big)_i \big(h_1'^{-1}\big)_{\pi(i)}
\big(h_2\big)_{\pi(i)} \big(h_2'^{-1}\big)_{\rho(i)}
\big(h_3\big)_{\rho(i)} = 1.
$$

By the simultaneous triple product property, we find that $\pi(i)
= \rho(i) = i$.  Thus, $\pi = \rho = 1$, which together with
\eqref{equation:stppperm} implies $\pi_i=\pi'_i$ for all $i$.
Finally, we have
$
h_1 h_1'^{-1} h_2 h_2'^{-1} h_3 h_3'^{-1} = 1,
$
which implies $h_1=h_1'$, $h_2=h_2'$, and $h_3=h_3'$ because each
triple $A_i,B_i,C_i$ satisfies the triple product property.
\end{proof}

\remove{ As a first step towards proving
Theorem~\ref{theorem:asi}, we prove a weaker bound, with the
geometric mean replacing the arithmetic mean:

\begin{lemma}
\label{lemma:wreath-bound} If $H$ is a finite group with character
degrees $\{d_k\}$ and $n$ triples of subsets $A_i, B_i, C_i
\subseteq H$ satisfying the simultaneous triple product property,
then
\[n\left(\prod_i(|A_i||B_i||C_i|)^{\omega/3}\right )^{1/n} \le \sum_k
d_k^\omega.\]
\end{lemma}

\begin{proof}
The sizes of the three subsets of $G$ in Theorem
\ref{theorem:STPP2TPP} are $n!\prod_i|A_i|$, $n!\prod_i|B_i|$, and
$n!\prod_i|C_i|$, respectively. Applying
Theorem~\ref{theorem:bound} we get the inequality
$$
\left ((n!)^3 \prod_i|A_i||B_i||C_i|\right )^{\omega/3}\le \sum_j
c_j^\omega.
$$
By Lemma \ref{lemma:wreath-char-degrees} the right-hand-side is at
most $(n!)^{\omega - 1}\left (\sum_k d_k^{\omega}\right )^n$, and
then dividing both sides by $(n!)^\omega$ yields
\[\left
(\prod_i|A_i||B_i||C_i|\right )^{\omega/3} \le (n!)^{-1}\left (
\sum_k d_k^\omega\right )^n.\]

This inequality is slightly weaker than the desired inequality,
but that is easy to fix by taking direct powers of $H$ via
Lemma~\ref{lemma:directprod}. Replacing $H$ with $H^t$ (and $n$
with $n^t$) yields
\[
\left (\prod_i|A_i||B_i||C_i|\right )^{tn^{t-1}\omega/3} \le
(n^t!)^{-1}\left ( \sum_k d_k^\omega\right )^{tn^t}.
\]
Taking $tn^t$-th roots and letting $t \to \infty$ gives the
claimed inequality.
\end{proof}
}

The proof of Theorem~\ref{theorem:asi} follows the same outline as
the proof of Theorem~\ref{theorem:SDPP-bound}:

\begin{proof}[Proof of Theorem~\ref{theorem:asi}]
Let $A_i', B_i', C_i'$ be the $N$-fold direct product (via Lemma
\ref{lemma:directprod}) of the triples $A_i, B_i, C_i$ realizing
$\langle a_1,b_1,c_1\rangle,\dots,\langle a_n,b_n,c_n \rangle$,
and let $\mu$ be an arbitrary $n$-vector of nonnegative integers
for which $\sum_{i=1}^n \mu_i = N$. Among the triples $A_i', B_i',
C_i'$ are $M = \binom{N}{\mu}$ triples for which $
|A_i'||B_i'||C_i'| = \prod_{i = 1}^n (a_ib_ic_i)^{\mu_i};$ call
this quantity $L$.

Using these triples in Theorem \ref{theorem:STPP2TPP}, we obtain
subsets $H_1, H_2, H_3$ with $|H_1||H_2||H_3| = (M!)^3L^M$.
Applying Theorem~\ref{theorem:bound} and Lemma
\ref{lemma:wreath-char-degrees} we get $((M!)^3L^M)^{\omega/3} \le
(M!)^{\omega - 1}\left ( \sum_k d_k^{\omega}\right )^{NM}$. Taking
$M$-th roots and letting $N \rightarrow \infty$ yields
$$
\binom{N}{\mu}\left (\prod_{i = 1}^n (a_ib_ic_i)^{\mu_i}\right
)^{\omega/3} \le \left (\sum_k d_k^\omega \right )^{N}.
$$
Finally, we apply Lemma \ref{lemma:geom->arith} with $s_i
=(a_ib_ic_i)^{\omega/3}$ and $C = \sum_k d_k^\omega$ to obtain the
stated inequality.
\end{proof}

\remove{
\begin{proof}[Proof of Theorem~\ref{theorem:asi}]
Let $A_i', B_i', C_i'$ be the $N$-fold direct product of the
triples $A_i, B_i, C_i$ (via Lemma \ref{lemma:directprod}), and
let $\mu$ be an arbitrary $n$-vector of nonnegative integers for
which $\sum_{i=1}^n \mu_i = N$. Among the triples $A_i', B_i',
C_i'$ are $\binom{N}{\mu}$ triples for which
$$
|A_i'||B_i'||C_i'| = \prod_{i = 1}^n (|A_i||B_i||C_i|)^{\mu_i}.
$$
Applying Lemma \ref{lemma:wreath-bound} to these triples gives
$$
\binom{N}{\mu}\prod_{i=1}^n(|A_i||B_i||C_i|)^{\mu_i\omega/3} \le
\left (\sum_k d_k^\omega\right )^N.
$$
Applying Lemma~\ref{lemma:geom->arith} with $s_i
=(|A_i||B_i||C_i|)^{\omega/3}$ and $C = \sum_k d_k^\omega$ yields
the desired bound.
\end{proof}
}

\section*{Acknowledgements}

We are grateful to Michael Aschbacher, Noam Elkies, William
Kantor, L\'aszl\'o Lov\'asz, Amin Shokrollahi, G\'abor Simonyi,
and David Vogan for helpful discussions.


\begin{thebibliography}{10}\setlength{\itemsep}{-1ex}\small

\bibitem{BCS}
P.~{B\"urgisser}, M.~Clausen, and M.~A. Shokrollahi.
\newblock {\em Algebraic Complexity Theory}, volume 315 of {\em Grundlehren der
  mathematischen Wissenschaften}.
\newblock Springer-Verlag, 1997.

\bibitem{CU}
H.~Cohn and C.~Umans.
\newblock A Group-theoretic Approach to Fast Matrix Multiplication.
\newblock Proceedings of the 44th Annual Symposium on Foundations of Computer Science,
11--14 October 2003, Cambridge, MA, IEEE Computer Society,
pp.~438--449, \texttt{arXiv:math.GR/0307321}.

\bibitem{CW}
D.~Coppersmith and S.~Winograd.
\newblock Matrix multiplication via arithmetic progressions.
\newblock {\em J. Symbolic Computation}, 9:251--280, 1990.


\bibitem{H}
B.~Huppert.
\newblock {\em Character Theory of Finite Groups}.
\newblock Number~25 in de Gruyter Expositions in Mathematics. Walter de
  Gruyter, Berlin, 1998.

\bibitem{JL}
G.~James and M.~Liebeck.
\newblock {\em Representations and Characters of Groups}.
\newblock Cambridge University Press, Cambridge, second edition, 2001.

\bibitem{SS}
R.~Salem and D.~C.~Spencer.
\newblock On sets of integers which contain no three terms in arithmetical progression.
\newblock {\em Proc.\ Nat.\ Acad.\ Sci.\ USA}, 28:561--563, 1942.

\bibitem{Simonyi}
G.~Simonyi.
\newblock Perfect graphs and graph entropy. An updated
survey.
\newblock P\textit{erfect graphs\/}, 293--328, Wiley-Intersci.\
Ser.\ Discrete Math.\ Optim., Wiley, Chichester, 2001.

\bibitem{S}
V.~Strassen.
\newblock Gaussian elimination is not optimal.
\newblock {\em Numerische Mathematik}, 13:354--356, 1969.

\bibitem{S87}
V.~Strassen.
\newblock Relative bilinear complexity and matrix multiplication.
\newblock {\em J.\ Reine Angew.\ Math.\/}, 375/376:406--443, 1987.

\end{thebibliography}
\end{document}